%% file: main.tex
\documentclass[12pt]{article}
\usepackage{CJKutf8}
\usepackage{amsmath,amsthm,amssymb,amsfonts}
\usepackage{color,colordvi}
\usepackage{soul}
\usepackage{fullpage}
\usepackage{amsmath,amsfonts,amsthm,mathtools,bm}
\usepackage{graphicx}
\usepackage[colorlinks=true, allcolors=blue]{hyperref}
\usepackage[capitalise]{cleveref}
\usepackage{standalone}
\usepackage{indentfirst}
\usepackage[affil-it]{authblk}
\usepackage{diagbox}
\usepackage{changes}
\usepackage{enumitem}
\usepackage{bbm}
\usepackage{orcidlink}

\newtheorem{theorem}{Theorem}[section]
\newtheorem{lemma}[theorem]{Lemma}

\newtheorem{proposition}[theorem]{Proposition}

\newtheorem{problem}[theorem]{Problem}
\newtheorem{conjecture}[theorem]{Conjecture}

\newtheorem{definition}[theorem]{Definition}

\newtheorem*{claim*}{Claim}

\newcommand{\textotherwise}{\text{otherwise}}

\newcommand{\CC}{\mathbb{C}}
\newcommand{\RR}{\mathbb{R}}

\newcommand{\QQ}{\mathbb{Q}}
\newcommand{\ZZ}{\mathbb{Z}}

\newcommand{\cP}{\mathcal{P}}

\newcommand{\fA}{\mathfrak{A}}
\newcommand{\fB}{\mathfrak{B}}

\newcommand{\finer}{\trianglelefteq}
\newcommand{\coarser}{\trianglerighteq}

\DeclareMathOperator{\Aut}{Aut}
\DeclareMathOperator{\cl}{cl}
\DeclareMathOperator{\tp}{atp}
\DeclareMathOperator{\Tr}{Tr}
\DeclareMathOperator{\diag}{diag}
\DeclareMathOperator{\offdiag}{offdiag}
\DeclareMathOperator{\Hom}{Hom}

\DeclarePairedDelimiter{\card}{\lvert}{\rvert}
\DeclarePairedDelimiter{\set}{\lbrace}{\rbrace}

\DeclarePairedDelimiter{\paren}{\lparen}{\rparen}

\DeclarePairedDelimiter{\lrangle}{\langle}{\rangle}

\DeclarePairedDelimiter{\ceil}{\lceil}{\rceil}

\title{Graph isomorphism and multivariate graph spectrum}

\author[a]{Wei Wang(\begin{CJK*}{UTF8}{gbsn}王卫\end{CJK*})}
\affil[a]{School of Mathematics and Statistics, Xi'an Jiaotong University, Xi'an, China.
\href{mailto:wang\_weiw@xjtu.edu.cn}{wang\_weiw@xjtu.edu.cn}}
\author[b]{Da Zhao(\begin{CJK*}{UTF8}{gbsn}赵达\end{CJK*})\orcidlink{0000-0002-9582-0778}}
\affil[b]{School of Mathematics, East China University of Science and Technology, Shanghai, China. \href{mailto:zhaoda@ecust.edu.cn}{zhaoda@ecust.edu.cn}}

\date{}

\begin{document}
\maketitle

\begin{abstract} 
    We provide a criterion to show that a graph is identified by its multivariate graph spectrum. 
    Haemers conjectured that almost all graphs are identified by their spectra. 
    Our approach suggests that almost all graphs are identified by their generalized block Laplacian spectra.  
\end{abstract}



\section{Introduction}

\begin{figure}[!htbp]
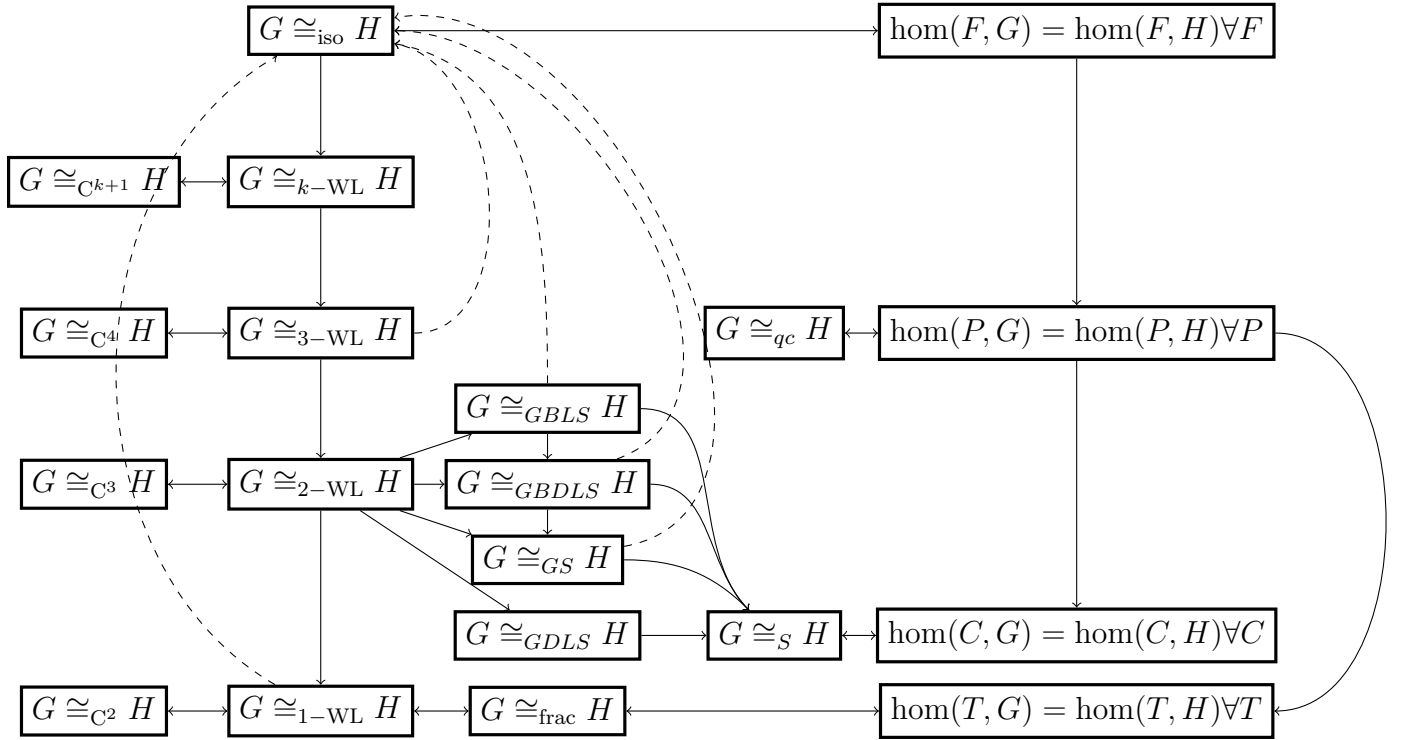

    \includestandalone{relation}
    \caption{Relation among equivalence classes of graphs}\label{fig:relation}
\end{figure}

The graph isomorphism problem asks whether two given graphs are isomorphic or not. 
Complexity theoretic results suggest that the isomorphism problem is not NP-complete~\cite{boppana1987DoesCoNPHave, schoning1988GraphIsomorphismLow}. 
On the other hand, no polynomial time algorithm is known so far. 
The color refinement algorithm is an intuitive algorithm that tries to tell whether two given graphs are isomorphic or not through iteratively coloring the vertices via the information on neighborhoods. 
Babai--Ku\v{c}era~\cite{babai1979CanonicalLabellingGraphs} and Babai--Erd\H{o}s--Selkow~\cite{babai1980RandomGraphIsomorphism} show that the color refinement algorithm succeeds almost surely in telling non-isomorphic graphs apart.
An improvement of color refinement algorithm is the $k$-dimensional Weisfeiler--Leman ($k$-WL) algorithm whose time complexity is bounded by $\mathcal{O}(n^{k+1} \log n)$~\cite{immerman2019} and the iteration number is bounded by $\mathcal{O}(k n^{k-1} \log n)$~\cite{grohe2025IterationNumberWeisfeilerleman} for $k \geq 2$. 
In 1992, Cai--F\"{u}rer--Immerman found a close connection between the WL algorithm and definability in certain finite variable logics studied in descriptive complexity theory. 
They show that for every $k$ there is a pair of non-isomorphic graphs of size linear in $k$ and valency $3$ such that they cannot be distinguished by the $k$-WL algorithm~\cite{cai1992OptimalLowerBound}. 
Though the power of the $k$-WL algorithm is limited, the algorithm with $k$ constant is guaranteed to work when restricted to certain graph families, such as planar or graphs with bounded genus~\cite{grohe1998FixedpointLogicsPlanar,grohe2000IsomorphismTestingEmbeddable}. 
Distinct from the WL-algorithm, Hopcroft--Wong~\cite{hopcroft1974LinearTimeAlgorithm} constructs a linear-time isomorphism algorithm for planar graphs.
Tinhofer~\cite{tinhofer1991NoteCompactGraphs} shows that two graphs are indistinguishable by the $1$-WL algorithm if and only if they are fractionally isomorphic. 
Lov\'{a}sz~\cite{lovasz1967OperationsStructures} proves that either graph homomorphism counter $\hom(G, \cdot)$ and $\hom(\cdot, G)$ identifies graph $G$. 
It is folklore that two graphs are cospectral if and only if the numbers of graph homomorphisms from cycle graphs are identical. 
Man\v{c}inska--Roberson~\cite{mancinskaQuantumIsomorphismEquivalent2019} remarkably shows that two graphs are quantum isomorphic if and only if their homomorphism counts from planar graphs are identical. 
It is known that two graphs are fractionally isomorphic if and only if their homomorphism counts from trees are identical. 
Haemers conjectures that almost all graphs are identified by their spectra. 

\subsection{Notions and notations}
To state our main result, we fix some notions and notations. A \emph{graph} $G$ is a pair $(V, E)$, where $V$ is the vertex set, and $E \subseteq \binom{V}{2}$ is the edge set. 
We denote an edge $(u, v) \in E$ by $uv$.
Note that here the graph is \emph{simple}, namely it has no loops, no multiedges, and the edges are undirected. 
The degree $\deg(v)$ of a vertex $v \in V$ is the number of edges incident to $v$.
The \emph{adjacency matrix} of a graph $G$ is a zero-one square matrix $A_G$ over the field $\CC$ given by
\begin{align}
    A_G(u, v) = 
    \begin{cases}
        1, & uv \in E, \\
        0, & \textotherwise. 
    \end{cases}
\end{align}
Let $G = (V_G, E_G)$ and $H = (V_H, E_H)$ be two graphs.
A \emph{graph homomorphism} from $G$ to $H$ is a map $f : V_G \to V_H$ such that $uv \in E_G$ implies $f(u)f(v) \in E_H$. 
The \emph{homomorphism count} $\hom(G, H)$ is the cardinality of $\Hom(G, H)$, where $\Hom(G, H)$ is the set of all graph homomorphisms from $G$ to $H$. 
A \emph{graph isomorphism} from $G$ to $H$ is a map $f : V_G \to V_H$ such that $uv \in E_G$ if and only if $f(u)f(v) \in E_H$. 
We write $G \cong_{\text{iso}} H$ if $G$ and $H$ are isomorphic. 
Two graphs $G$ and $H$ are fractionally isomorphic, denoted by $G \cong_{\text{frac}} H$, if there exists a \emph{doubly stochastic matrix} $X$ such that $A_G X = X A_H$, where a doubly stochastic matrix is a nonnegative square matrix such that the sum of each row or column is $1$. 
In particular, a \emph{permutation matrix} is a doubly stochastic matrix. 

A matrix is called \emph{integral} if all of its entries are integers. 
A square integral matrix $M$ is called \emph{unimodular} if $\det M=\pm1$. Note that for a unimodular matrix $M$, its inverse $M^{-1}$ is always  an integral matrix. 
Given a rectangular $n \times m$ integral matrix $M$, there exists a decomposition $M = U \Sigma V$, called \emph{Smith decomposition}, such that the followings hold. 
\begin{enumerate}
    \item The matrix $U$ is an $n \times n$ unimodular matrix;
    \item The matrix $V$ is an $m \times m$ unimodular matrix;
    \item The matrix $\Sigma$ is a $n \times m$ diagonal matrix such that the diagonal entries satisfy $d_1 \mid d_2 \mid \cdots \mid d_{\min(n,m)}$.
\end{enumerate}
We call $\Sigma$ the \emph{Smith normal form} of $M$ and $d_{\min(n,m)}$ the \emph{last invariant factor} of $M$. 

Given a rectangular matrix $M$ over $\QQ$, the \emph{level} of $M$, denoted by $\ell(M)$, is the smallest positive integer $\ell$ such that $\ell M$ is an integral matrix. 

Let $\lambda_1, \lambda_2, \ldots, \lambda_n$ be the eigenvalues of a square matrix $M$ of order $n$, the \emph{discriminant} of $M$ is defined by $\Delta_M = \prod_{i<j} (\lambda_i - \lambda_j)^2$. 
Note that if $A$ is an integral matrix, then $\Delta_A$ is always a non-negative integer.

Throughout, we shall denote by $I$ and $J$  the identity matrix and the all-one matrix respectively. 

\subsection{Main results} 

The main theorem of this paper is as follows. 

\begin{theorem}\label{thm:main}
    Let $G$ and $H$ be two graphs of order $n$. 
    Suppose the $2$-WL algorithm cannot distinguish $G$ and $H$, namely $G \cong_{2\text{-WL}} H$. 
    Let $A$ and $B$ be the adjacency matrices of $G$ and $H$ respectively. 
    Suppose the partition of vertices of $G$ by degree is given by $V = \bigsqcup_{i=1}^p V_i$, and the characteristic vector corresponding to $V_i$ is $e_i$. 
    If the last invariant factor $d_n(\widetilde{W}_A)$ in Smith normal form of $\widetilde{W}_A = [e_1, A e_1, \ldots, A^{n-1} e_1, e_2, A e_2, \ldots, A^{n-1} e_2, \ldots, e_p, A e_p, \ldots, A^{n-1} e_p]$ is odd and for every odd prime $q$ dividing $d_n(\widetilde{W}_A)$, it holds that $q^2$ does not divide the discriminant of $A$, then $G \cong_{\text{iso}} H$. 
    In particular, if the last invariant factor of $\widetilde{W}_A$ is $1$, then $G \cong_{\text{iso}} H$.
    The above claims hold for $B$ as well. 
\end{theorem}

Let $G$ be a graph. 
The (adjacency) \emph{spectrum} of $G$ is the multiset of eigenvalues of $A_G$.
We say two graphs $G$ and $H$ are \emph{cospectral}, denoted by $G \cong_{S} H$, if their spectra are identical. 
It is well-known that $2$-WL equivalence implies cospectral.  
In fact, we can consider more general spectrum. 
The \emph{complementary graph} $\overline{G}$ of $G = (V, E)$ is a graph $(V, \overline{E})$ such that $\overline{E} \coloneqq \binom{V}{2} \setminus E$. 
The \emph{generalized spectrum} of $G$ is the collection of the spectrum of $G$ as well as the spectrum of $\overline{G}$. 
We say $G$ and $H$ are \emph{generalized cospectral}, denoted by $G \cong_{GS} H$, if their generalized spectra are identical. 

Let $\bm{A} = (A_1, A_2, \ldots, A_k)$ be a $k$-tuple of $n \times n$ integral matrices. 
Let $\bm{s} = (s_1, s_2, \ldots, s_k)$ be a $k$-tuple of complex variables. 
We define
\begin{align}
    W_{\bm{A}}(\bm{s}) = W_{A_1, A_2, \ldots, A_k}(s_1, s_2, \ldots, s_k) = \sum_{i=1}^k s_i A_i.
\end{align}
Denote by $\phi_{\bm{A}}(\bm{s}; t) = \det(t I - W_{\bm{A}}(\bm{s}))$ the \emph{characteristic polynomial} of $W_{\bm{A}}(\bm{s})$. 
It is clear that $G$ and $H$ are cospectral if and only if $\phi_{\bm{A}}(\bm{s}; t) = \phi_{\bm{B}}(\bm{s}; t)$, where $\bm{A} = (A_G)$ and $\bm{B} = (A_H)$.
In fact, we also have that $G$ and $H$ are generalized cospectral if and only if $\phi_{\bm{A}}(\bm{s}; t) = \phi_{\bm{B}}(\bm{s}; t)$, where $\bm{A} = (A_G, J)$ and $\bm{B} = (A_H, J)$~\cite{johnson1980NoteCospectralGraphs}.
In the following we will generalize the above fact. 

Let $G$ and $H$ be two graphs on $V$ sharing the same degree sequences. 
Without loss of generality, suppose $\deg(v) = \deg_{G}(v) = \deg_{H}(v)$ for all $v \in V$. 
Therefore, we have the degree decomposition of vertices $V = \bigsqcup_{i=1}^p V_i$.
Let $A$ and $B$ be the adjacency matrices of $G$ and $H$ respectively. 
Let $e_i \in \RR^V$, $i = 1,2, \ldots, p$ be zero-one vectors such that $e_i(v) = 1$ if and only if $v \in V_i$. 
Let $J_{i,j} = e_i e_j^\top$ for $i,j = 1, 2, \ldots, p$.
Let $D_i(u, v) = 
\begin{cases}
    1, & u = v \in V_i, \\
    0, & \textotherwise 
\end{cases}$.
We say two graphs $G$ and $H$ share the \emph{generalized diagonal Laplacian spectrum}, denoted by $G \cong_{GDLS} H$, if $\phi_{\bm{A}'}(\bm{s}, t) = \phi_{\bm{B}'}(\bm{s}, t)$, where $\bm{A}' = (A_G, D_1, D_2, \ldots, D_p)$ and $\bm{B}' = (A_H, D_1, D_2, \ldots, D_p)$.
We say two graphs $G$ and $H$ share the \emph{generalized block diagonal Laplacian spectrum}, denoted by $G \cong_{GBDLS} H$, if $\phi_{\bm{A}''}(\bm{s}, t) = \phi_{\bm{B}''}(\bm{s}, t)$, where $\bm{A}'' = (A_G, J_{1,1}, J_{2,2}, \ldots, J_{p,p})$ and $\bm{B}'' = (A_H, J_{1,1}, J_{2,2}, \ldots, J_{p,p})$.
We say two graphs $G$ and $H$ share the \emph{generalized block Laplacian spectrum}, denoted by $G \cong_{GBLS} H$, if $\phi_{\bm{A}'''}(\bm{s}, t) = \phi_{\bm{B}'''}(\bm{s}, t)$, where $\bm{A}''' = (A_G, J_{1,1}, J_{1,2}, \ldots, J_{p,p})$ and $\bm{B}''' = (A_H, J_{1,1}, J_{1,2}, \ldots, J_{p,p})$. 

\cref{thm:main} enables us to show the isomorphism among some graphs without specifying the map between the vertices.  
In fact, the last invariant factor of $\widetilde{W}_A$ suggests potential graphs sharing the generalized block Laplacian spectrum. 
All such graphs can be found, with minor modifications, by the algorithm in~\cite{wang2024HaemersConjectureAlgorithmic} as long as the last invariant factor is not zero.
\cref{thm:main} is a corollary of the following three theorems. 

\begin{theorem}\label{thm:2WLtoGBLS}
    Let $G$ and $H$ be two graphs of order $n$. 
    Suppose the $2$-WL algorithm cannot distinguish $G$ and $H$, namely $G \cong_{2\text{-WL}} H$. 
    Then $G$ and $H$ share the generalized block Laplacian spectrum, namely $G \cong_{GBLS} H$. 
\end{theorem}

\begin{theorem}\label{thm:GBLStodn}
    Let $G$ and $H$ be two graphs of order $n$.
    Suppose $G$ and $H$ share the generalized block Laplacian spectrum, namely $G \cong_{GBLS} H$. 
    Let $V = \bigsqcup_{i=1}^p V_i$ be the partition of vertices of $G$ by degree, and let $e_i$ be the characteristic vector corresponding to $V_i$
    Then there exists a fixed orthogonal matrix $Q$ such that $Q^\top A Q = B$, and $Q^\top e_i = e_i$. 
    In particular, if $\widetilde{W}_A$ is of full row rank, then $Q$ is unique and rational. 
    Moreover, $d_n(\widetilde{W}_A) Q$ is an integral matrix, where $d_n(X)$ is the last invariant factor of $X$. 
\end{theorem}

\begin{theorem}[{\cite[Lemma 4.5]{wang2016SquarefreeDiscriminantsMatrices}}]
    Let $Q$ be a rational orthogonal matrix. Let $A$ be a symmetric integral matrix such that $Q^\top A Q$ is a symmetric integral matrix. 
    Suppose $p$ is an odd prime, $p \mid \ell(Q)$, and $p \mid \Delta_A$, where $\Delta_A$ is the discriminant of $A$.  
    Then $p^2 \mid \Delta_A$. 
\end{theorem}

We summarize the relation among some equivalence classes of graphs in~\cref{fig:relation}, where the symbols $F$, $P$, $C$, and $T$ are finite graphs, planar graphs, cycle graphs, and trees respectively. 
The straight lines indicate implications. 
The dashed lines indicate conditional implications.

\section{Multivariate graph spectrum}\label{sec:preliminary}

\begin{theorem}\label{thm:GBLS}
    Let $A_1$ and $B_1$ be two real symmetric matrices of order $n$. 
    Let $e_i$, $i = 1,2, \ldots, p$ be zero-one vectors of length $n$ such that the positions of the ones are disjoint. 
    Let $J_{i,j} = e_i e_j^\top$ for $i,j = 1, 2, \ldots, p$ and $k = 1 + p^2$.
    Let $\bm{A} = (A_1, J_{1,1}, J_{1,2}, \ldots, J_{p,p})$ and $\bm{B} = (B_1, J_{1,1}, J_{1,2}, \ldots, J_{p,p})$. 
    Then the matrices $W_{\bm{A}}(\bm{s})$ and $W_{\bm{B}}(\bm{s})$ have identically the same characteristic polynomials if and only if they are similar via a fixed orthogonal matrix $Q$, independent of $\bm{s}$, such that $Q^\top A_1 Q = B_1$ and $Q^\top e_i = e_i$ for $i = 1,2, \ldots, p$.
    Moreover, $Q^\top \widetilde{W}_{A_1} = \widetilde{W}_{B_1}$, where
    \begin{align}
        \widetilde{W}_{A} = [e_1, Ae_1, \ldots, A^{n-1} e_1, e_2, Ae_2, \ldots, A^{n-1} e_2, \ldots, e_p, A e_p, \ldots, A^{n-1} e_p].
    \end{align}
\end{theorem}

\begin{proof}
    On one hand, suppose there exists a fixed orthogonal matrix $Q$ such that $Q^\top A_1 Q = B_1$ and $Q^\top e_i = e_i$ for $i = 1,2, \ldots, p$. 
    Then $Q^\top J_{i,j} = J_{i,j} = Q J_{i,j}$ for $1 \leq i, j \leq p$. 
    Take transpose on both sides and we get $J_{i,j} Q = J_{i,j}$ for $1 \leq i, j \leq p$. 
    Hence, $\det(tI - W_{\bm{A}}(s)) = \det(Q^\top(tI - W_{\bm{A}}(s))Q) = \det(tI - Q^\top W_{\bm{A}}(s) Q) = \det(tI - W_{\bm{B}}(s))$. 
    Meanwhile, for $m = 0,1, \ldots, n-1$ we have 
    \begin{align}
        Q^\top A_1^m e_i = B_1^m Q^\top e_i = B_1^m e_i.
    \end{align}
    In other words $Q^\top \widetilde{W}_{A_1} = \widetilde{W}_{B_1}$.

    On the other hand, suppose $\det(tI - W_{\bm{A}}(s)) = \det(tI - W_{\bm{B}}(s))$.
    Note that there exists an orthogonal matrix $O$ such that
    \begin{align}
        O^\top J_{i,j} O = n_{i,j} E_{i,j},
    \end{align}
    where $E_{i,j}$ is the matrix whose only nonzero entry is a $1$ in the $(i,j)$-entry.
    Replace $\hat{A} = O^\top A O$ for $A \in \set{A_1, B_1, J_{1,1}, J_{1,2}, \ldots, J_{p,p}}$. 
    For $\hat{\bm{A}} = (\hat{A}_1, \hat{J}_{1,1}, \ldots, \hat{J}_{p,p})$ and $\hat{\bm{B}} = (\hat{B}_1, \hat{J}_{1,1}, \ldots, \hat{B}_{p,p})$, we have 
    \begin{align}\label{eqn:same_det}
        \det(t I - x \hat{A}_1 - \sum_{i,j=1}^p s_{i,j} n_{i,j} E_{i,j}) = \det(t I - x \hat{B}_1 - \sum_{i,j=1}^p s_{i,j} n_{i,j} E_{i,j})
    \end{align}
    for all $x, s_{1,1}, s_{1,2}, \ldots, s_{p,p} \in \CC$. 
    We regard~\cref{eqn:same_det} as a polynomial identity in $s_{i,j}, i,j = 1,2, \ldots, p$. 
    In the full expansion of $\det(t I - x \hat{A}_1 - \sum_{i,j=1}^p s_{i,j} n_{i,j} E_{i,j})$ by multilinearity, the coefficient for $s_{1,1} s_{2,2} \cdots s_{p,p}$ is
    \begin{align}
        &(-1)^p n_{1,1} n_{2,2} \cdots n_{p,p} \det(tI - x \hat{A}_1 - \sum_{i,j=1}^p s_{i,j} n_{i,j} E_{i,j})_{[F], [F]} \\
        =& (-1)^p n_{1,1} n_{2,2} \cdots n_{p,p} \det(tI - x \hat{A}_1)_{[F], [F]},
    \end{align}
    where $F$ denotes the set $\set{1,2, \ldots, p}$ and $X_{[K], [L]}$ denotes the submatrix of $X$ obtained by deleting the rows corresponding to $K$ and the columns corresponding to $L$. 
    Similarly, in the full expansion of $\det(t I - x \hat{B}_1 - \sum_{i,j=1}^p s_{i,j} n_{i,j} E_{i,j})$, the coefficient of $s_{1,1} s_{2,2} \cdots s_{p,p}$ is $(-1)^p n_{1,1} n_{2,2} \cdots n_{p,p} \det(tI - x \hat{B}_1)_{[F], [F]}$.
    Therefore, $\det(tI - x \hat{A}_1)_{[F],[F]} = \det(tI - x \hat{B}_1)_{[F], [F]}$.
    It follows that $(\hat{A}_1)_{[F], [F]}$ and $(\hat{B}_1)_{[F], [F]}$ are orthogonally similar. 
    We may assume that after suitable similarities of the form $P = \begin{bmatrix}
        I_{p} & 0 \\
        0 & U
    \end{bmatrix}$ it holds
    \begin{align}
        \hat{A}_1 = 
        \begin{bmatrix}
            a_{1,1} & \cdots & a_{1,p} & \alpha_{p+1,1}^\top & \cdots & \alpha_{q, 1}^\top \\
            \vdots & \ddots & \vdots & \vdots & \cdots & \vdots \\
            a_{p,1} & \cdots & a_{p,p} & \alpha_{p+1, p}^\top & \cdots & \alpha_{q, k}^\top \\
            \alpha_{p+1,1} & \cdots & \alpha_{p+1, p} & \lambda_{p+1} I_{m_{p+1}} & 0 & 0 \\
            \vdots & \vdots & \vdots & 0 & \ddots & 0 \\
            \alpha_{q,1} & \cdots & \alpha_{q,p} & 0 & 0 & \lambda_q I_{m_q}
        \end{bmatrix}
    \end{align}
    and
    \begin{align}
        \hat{B}_1 = 
        \begin{bmatrix}
            b_{1,1} & \cdots & b_{1,p} & \beta_{p+1,1}^\top & \cdots & \beta_{q, 1}^\top \\
            \vdots & \ddots & \vdots & \vdots & \cdots & \vdots \\
            b_{p,1} & \cdots & b_{p,p} & \beta_{p+1, p}^\top & \cdots & \beta_{q, k}^\top \\
            \beta_{p+1,1} & \cdots & \beta_{p+1, p} & \lambda_{p+1} I_{m_{p+1}} & 0 & 0 \\
            \vdots & \vdots & \vdots & 0 & \ddots & 0 \\
            \beta_{q,1} & \cdots & \beta_{q,p} & 0 & 0 & \lambda_q I_{m_q}
        \end{bmatrix}, 
    \end{align}
    where $\lambda_{p+1}, \ldots, \lambda_q$ are the distinct eigenvalues of $(\hat{A}_1)_{[F], [F]}$ (and also of $(\hat{B}_1)_{[F], [F]}$). 
    Note that $P^\top E_{i,j} P = E_{i,j}$ for $i,j = 1,2, \ldots, p$. 

    In the full expansion of $\det(t I - x \hat{A}_1 - \sum_{i,j=1}^p s_{i,j} n_{i,j} E_{i,j})$, the coefficient for $s_{2,2} \cdots s_{p,p}$ is $(-1)^p n_{2,2} \cdots n_{p,p} \det(tI - x \hat{A}_1 - \sum_{i,j=1}^p s_{i,j} n_{i,j} E_{i,j})_{[\set{2,3, \ldots, p}], [\set{2,3,\ldots, p}]}$. 
    So 
    \begin{align}
        &\det(tI - x \hat{A}_1 - \sum_{i,j=1}^p s_{i,j} n_{i,j} E_{i,j})_{[\set{2,3, \ldots, p}], [\set{2,3,\ldots, p}]} \\
        =& \det(tI - x \hat{B}_1 - \sum_{i,j=1}^p s_{i,j} n_{i,j} E_{i,j})_{[\set{2,3, \ldots, p}], [\set{2,3,\ldots, p}]}.
    \end{align}

    For $i,\ell \in F$, we define
    \begin{align}
        \fA_{i,\ell} = \delta_{i,\ell} t - x a_{i,\ell} - s_{i, \ell} n_{i, \ell} + \sum_{j=p+1}^q \frac{x^2 \lrangle{\alpha_{j,i}, \alpha_{j,\ell}}}{t - x \lambda_j},
    \end{align}
    and
    \begin{align}
        \fB_{i,\ell} = \delta_{i,\ell} t - x b_{i,\ell} - s_{i, \ell} n_{i, \ell} + \sum_{j=p+1}^{q} \frac{x^2 \lrangle{\beta_{j,i}, \beta_{j,\ell}}}{t - x \lambda_j},
    \end{align}
    where $\delta_{i, \ell}$ is the Kronecker delta function, namely
    \begin{align}
        \delta_{i,\ell} = 
        \begin{cases}
            1, & i = \ell, \\
            0, & i \neq \ell.
        \end{cases}
    \end{align}
    By elementary row transformation, we get
    \begin{align}
        \det(t I - x \hat{A}_1  - \sum_{i,j=1}^p s_{i,j} n_{i,j} E_{i,j})_{[\set{2,3,\ldots, p}],[\set{2,3,\ldots, p}]} &= \fA_{1,1} \prod_{j=p+1}^q (t - x \lambda_j)^{m_j}
    \end{align}
    and
    \begin{align}
        \det(t I - x \hat{B}_1  - \sum_{i,j=1}^p s_{i,j} n_{i,j} E_{i,j})_{[\set{2,3,\ldots, p}],[\set{2,3,\ldots, p}]} &= \fB_{1,1} \prod_{j=p+1}^q (t - x \lambda_j)^{m_j}.
    \end{align}
    We conclude that $a_{1,1} = b_{1,1}$ and $\lrangle{\alpha_{j,1}, \alpha_{j,1}} = \lrangle{\beta_{j,1}, \beta_{j,1}}$ for $j = p, p+1, \ldots, q$. 
    Similarly, we have $a_{i,i} = b_{i,i}$ and $\lrangle{\alpha_{j,i}, \alpha_{j,i}} = \lrangle{\beta_{j,i}, \beta_{j,i}}$ for $i = 1, 2, \ldots, p$ and $j = p, p+1, \ldots, q$. 

    Now we consider the coefficient for $s_{3,3} \cdots s_{p,p}$ in the full expansion of $\det(t I - x \hat{A}_1 - \sum_{i,j=1}^p s_{i,j} n_{i,j} E_{i,j})$ as well as $\det(t I - x \hat{B}_1 - \sum_{i,j=1}^p s_{i,j} n_{i,j} E_{i,j})$, and get
    \begin{align}
        \det(t I - x \hat{A}_1  - \sum_{i,j=1}^p s_{i,j} n_{i,j} E_{i,j})_{[\set{3, 4, \ldots, p}],[\set{3, 4, \ldots, p}]} 
        = (\fA_{1,1} \fA_{2,2} - \fA_{1,2} \fA_{2,1}) \prod_{j=p+1}^q (t - x \lambda_j)^{m_j},
    \end{align}
    and
    \begin{align}
        \det(t I - x \hat{B}_1  - \sum_{i,j=1}^p s_{i,j} n_{i,j} E_{i,j})_{[\set{3, 4, \ldots, p}],[\set{3, 4, \ldots, p}]} 
        = (\fB_{1,1} \fB_{2,2} - \fB_{1,2} \fB_{2,1}) \prod_{j=p+1}^q (t - x \lambda_j)^{m_j}.
    \end{align}
    Therefore, 
    \begin{align}
        \fA_{1,1} \fA_{2,2} - \fA_{1,2} \fA_{2,1} = \fB_{1,1} \fB_{2,2} - \fB_{1,2} \fB_{2,1}.
    \end{align}
    Consider the polynomial term, and we get
    \begin{align}
        & \paren*{t - x a_{1,1} - s_{1,1} n_{1,1}} \paren*{t - x a_{2,2} - s_{2,2} n_{2,2}} 
        - \paren*{- x a_{1,2} - s_{1,2} n_{1,2}} \paren*{- x a_{2,1} - s_{2,1} n_{2,1}}  \\
        =
        & \paren*{t - x b_{1,1} - s_{1,1} n_{1,1}} \paren*{t - x b_{2,2} - s_{2,2} n_{2,2}} 
        - \paren*{- x b_{1,2} - s_{1,2} n_{1,2}} \paren*{- x b_{2,1} - s_{2,1} n_{2,1}}.
    \end{align}
    Hence, $a_{1,2} = b_{1,2}$ and $a_{2,1} = b_{2,1}$. 
    Similarly, we have $a_{i,j} = b_{i,j}$ for $1 \leq i, j \leq p$. 
    Consider the coefficient for $\frac{1}{t - x \lambda_j}$, and we get
    \begin{align}
        & \paren*{t - x a_{1,1} - s_{1,1} n_{1,1}}(x^2 \lrangle{\alpha_{j,2}, \alpha_{j,2}}) + \paren*{t - x a_{2,2} - s_{2,2} n_{2,2}}(x^2 \lrangle{\alpha_{j,1}, \alpha_{j,1}}) \\
        &- \paren*{- x a_{1,2} - s_{1,2} n_{1,2}}(x^2 \lrangle{\alpha_{j,2}, \alpha_{j,1}}) - \paren*{- x a_{2,1} - s_{2,1} n_{2,1}}(x^2 \lrangle{\alpha_{j,1}, \alpha_{j,2}}) \\
        =& \paren*{t - x b_{1,1} - s_{1,1} n_{1,1}}(x^2 \lrangle{\beta_{j,2}, \beta_{j,2}}) + \paren*{t - x b_{2,2} - s_{2,2} n_{2,2}}(x^2 \lrangle{\beta_{j,1}, \beta_{j,1}}) \\
        &- \paren*{- x b_{1,2} - s_{1,2} n_{1,2}}(x^2 \lrangle{\beta_{j,2}, \beta_{j,1}}) - \paren*{- x b_{2,1} - s_{2,1} n_{2,1}}(x^2 \lrangle{\beta_{j,1}, \beta_{j,2}}).
    \end{align}
    Hence, $\lrangle{\alpha_{j,1}, \alpha_{j,2}} = \lrangle{\beta_{j,1}, \beta_{j,2}}$, $\lrangle{\alpha_{j,2}, \alpha_{j,1}} = \lrangle{\beta_{j,2}, \beta_{j,1}}$ for $j = p+1, p+2, \ldots, q$. 
    Similarly, we have $\lrangle{\alpha_{j,i}, \alpha_{j, \ell}} = \lrangle{\beta_{j,i}, \beta_{j,\ell}}$, for $j = p+1, p+2, \ldots, q$, $1 \leq i \neq \ell \leq p$. 
    In other words, the inner products among the $p$ vectors $\alpha_{j, 1}, \alpha_{j, 2}, \ldots, \alpha_{j, p}$ are identical to those among the $p$ vectors $\beta_{j, 1}, \beta_{j, 2}, \ldots, \beta_{j, p}$ for fixed $j \in \set{p+1, p+2, \ldots, q}$. 
    So there exists an orthogonal transformation $V_j$ such that $V_j^\top \beta_{j, i} = \alpha_{j,i}$ for $i = 1,2, \ldots, p$. 
    Consider the orthogonal matrix $R = \begin{bmatrix}
        I_{p} & 0 & \cdots & 0\\
        0 & V_{p+1} & 0 & 0 \\
        0 & 0 & \ddots & 0 \\
        0 & 0 & 0 & V_{q} 
    \end{bmatrix}$. 
    Then we have $\hat{A}_1 = R^\top \hat{B}_1 R$ and $R^\top E_{i, j} R = E_{i,j}$ for $1 \leq i, j \leq p$. 
    Let $Q = O R O^\top$. 
    Then
    \begin{align}
        Q^\top A_1 Q = (ORO^\top)^\top O \hat{A}_1 O^\top ORO^\top = O R^\top \hat{A}_1 R O^\top = O\hat{B}_1 O^\top = B_1,
    \end{align}
    \begin{align}
        Q^\top J_{i, j} Q = (ORO^\top)^\top J_{i,j} ORO^\top = O R^\top (n_{i,j} E_{i,j}) R O^\top = O (n_{i,j} E_{i,j}) O^\top = J_{i,j},
    \end{align}
    and
    \begin{align}
        Q^\top J_{i,j} = O R^\top O^\top J_{i,j} = O R^\top (n_{i,j} E_{i,j}) O^\top = O (n_{i,j} E_{i,j}) O^\top = J_{i,j}.
    \end{align}
    Hence, $Q^\top W_{\bm{A}}(\bm{s}) Q = W_{\bm{B}}(\bm{s})$ and $Q^\top e_i = e_i$ for $i = 1, 2, \ldots, p$. 
\end{proof}

Note that in the above theorem $J_{i,j}$, $(1 \leq i,j \leq k)$ are all rank one matrices. 
We may consider other rank one matrices of the form $\xi \xi^\top$, where $\xi$ is a so-called graph-vector (See~\cite{qiu2023GeneralizedSpectralCharacterizations}). 

\begin{definition}
    Let $G = (V, E)$ be a graph. 
    A graph-vector of $G$ is a vector $\xi_G \in \RR^V$ such that it is invariant under the automorphisms of $G$, namely if $P$ is a permutation matrix such that $P^\top A_G P = A_G$, then $\xi_G = P^\top \xi_G$.
\end{definition}

We characterize the linear space spanned by all graph-vectors.

\begin{lemma}
    Let $G = (V, E)$ be a graph. 
    Then all the graph-vectors of $G$ form a linear subspace of $\RR^V$, whose dimension is equal to the number of orbits of $\Aut(G)$. 
    More precisely, the graph-vectors of $G$ are those which take constant value within an orbit. 
\end{lemma}

\begin{proof}
    We construct a basis of graph-vectors directly. 
    Let $(G,x_0)$ be a rooted graph where $x_0 \in V$. 
    Consider the vector $\xi_{G,x_0} \in \RR^V$ as follows. 
    \begin{align}
        \xi_{G, x_0}(x) = 
        \begin{cases}
            1, & (G, x) \cong (G, x_0) \\
            0, & (G, x) \not\cong (G, x_0).
        \end{cases}
    \end{align}
    It is clear that $\xi_{G, x_0}$ is invariant under the action of $\Aut(G)$.
    Therefore, $\xi_{G, x_0}$ is indeed a graph-vector. 
    Moreover, $\xi_{G, x_0} = \xi_{G, x_0'}$ if $x_0$ and $x_0'$ are in the same orbit of $\Aut(G)$. 

    Note that a graph-vector is invariant under the action of $\Aut(G)$, therefore it is a linear combination of $\xi_{G, x_0}$ with $x_0 \in V$. 
    We conclude that $\set{\xi_{G, x_0}, x_0 \in V}$ gives a basis of graph-vectors, whose size is equal to the number of orbits of $\Aut(G)$. 
\end{proof}

Next, we exhibit the restriction on the level of a rational matrix $Q$ by the identity $QX = Y$. 

The level of a rational matrix is invariant by the multiplication of a unimodular matrix. 

\begin{lemma}\label{lem:multiply_unimodular}
    Let $Q$ be a rational matrix and $U$ be a unimodular matrix. 
    Then $\ell(Q) = \ell(QU)$. 
\end{lemma}

\begin{proof}
    Suppose $R = QU$. 
    Then $\ell(Q)R = \ell(Q)Q U$ is an integral matrix. 
    Therefore, $\ell(R) \mid \ell(Q)$. 
    On the other hand, $Q = R U^{-1}$ and $U^{-1}$ is unimodular. 
    Therefore, $\ell(Q) \mid \ell(R)$. 
    Altogether, we get $\ell(Q) = \ell(QU)$.
\end{proof}

The identity $Q X= Y$ restricts the level of $Q$. 

\begin{lemma}\label{lem:divide_dn}
    Let $X$ and $Y$ be two $n \times m$ integral matrices ($n \leq m$) and let $Q$ be a square matrix of order $n$. 
    Suppose $Q X = Y$ and $X$ is of full row rank. 
    Then $Q$ is a unique rational matrix and  $\ell(Q) \mid d_n(X)$, where $d_n(X)$ is the last invariant factor of $X$.
\end{lemma}

\begin{proof}
    For each row $q^\top$ of $Q$ and the corresponding row $y^\top$ of $Y$, we have $q^\top X = y^\top$. 
    Since $X$ is of full rank and both $X$ and $Y$ are integral matrices, the row $q^\top$ is determined by $X$ and $Y$, and it is rational. 
    Hence, $Q$ is a unique rational matrix. 
    Consider the Smith decomposition of $X$, namely $X = U \Sigma V$, where both $U$ and $V$ are unimodular matrices. 
    Then $Q U \Sigma V = Y$. 
    Hence, $Q U \Sigma = Y V^{-1}$ is integral. 
    Note that $d_{n}(X) = d_{n}(\Sigma)$. 
    Hence, $\ell(Q) = \ell(QU) \mid d_n(\Sigma) = d_n(X)$.
\end{proof}

\section{The Weisfeiler-Leman algorithm}

Let $A : V \times V \to C$ be a square matrix indexed by $V$, whose entries take value in a set $C$. 
We define an equivalence relation on the set of $k$-tuples in $V^k$. 
Two $k$-tuples $(i_1, i_2, \ldots, i_k)$ and $(j_1, j_2, \ldots, j_k)$ are equivalent if
\begin{enumerate}
    \item $i_\ell = i_{\ell'}$ if and only if $j_\ell = j_\ell'$,
    \item $A(i_\ell, i_{\ell'}) = c$ if and only if $A(j_\ell, j_{\ell'}) = c$. 
\end{enumerate}

We define the atomic type $\tp(i_1, i_2, \ldots, i_k) = \tp_A(i_1, i_2, \ldots, i_k)$ of a $k$-tuple as its equivalence class. 
The information of the atomic type $\tp(i_1, i_2, \ldots, i_k)$ can be encoded in a $k \times k$ matrix $T$ with 
\begin{align}
    T_{p,q} = T(A, k)_{p,q} = 
    \begin{cases}
        (0, A(i_p, i_q)), & i_p = i_q, \\
        (1, A(i_p, i_q)), & i_p \neq i_q.
    \end{cases}
\end{align}
Let $S_1$ be the set of all different types of $k$-tuples, and it is the set of initial colors. 
We define the set $S$ of colors by
\begin{align}
    S = \bigcup_{k=1}^{\infty} S_k,
\end{align}
where elements of $S_{r+1}$ are finite multisets, which can be regarded as formal sums, of finite sequences whose elements are in $\cup_{k=1}^r S_k$. 
In fact, it is enough to work with as many colors as $k$-tuples by renaming the colors in each round. 

The $k$-WL algorithm works by iteratively assigning colors to $V^k$. 
That is to say, we define the functions $X_{A, k}^r : V^k \to S$. 
For $r=1$, we define
\begin{align}
    X_{A, k}^1 (i_1 \ldots i_k) = \tp_A(i_1 \ldots i_k).
\end{align}
Sequentially we have
\begin{align}\label{eq:iter_rplus1}
    X_{A, k}^{r+1}(i_1 \ldots i_k) = \paren*{X_{A, k}^r (i_1 \ldots i_k), \sum_{m \in V} (\tp_A(i_1 \ldots i_k m), S_{A, k}^r (i_1 \ldots i_k m))},
\end{align}
where $S_{A, k}^r (i_1 \ldots i_k m)$ is the sequence
\begin{align}
    (X_{A, k}^r (i_1 \ldots i_{k-1} m), \ldots, X_{A,k}^r (i_1 \ldots i_{\ell - 1} m i_{\ell + 1} \ldots i_k), \ldots, X_{A, k}^r (m i_2 \ldots i_k))
\end{align}
The summation in~\cref{eq:iter_rplus1} is regarded as a formal sum. 

Note that $X_{A, k}^r$ induces a partition $\cP(A, k, r)$ of the set $V^k$ of $k$-tuples for each $r \geq 1$. 
In particular, $\tp_A$, in other words, $X_{A, k}^1$, induces a partition, denoted by $\cP(A, k)$. 
The partition sequence stabilizes for sufficiently large $r$, denoted by $\cP(A, k, \infty)$. 
We define the closure $\cl(P)$ of $P = \cP(A, k)$ as $\cP(A, k, \infty)$. 
There is a partial order among the partitions of $V^k$. 
Let $P_1, P_2$ be two partitions of a set $X$. 
We say $P_1$ is finer than $P_2$, or $P_2$ is coarser than $P_1$, denoted by $P_1 \finer P_2$, if every part of $P_1$ is contained in a part of $P_2$. 
Note that
\begin{align}
    \cP(A, k) = \cP(A, k, 1) \coarser \cP(A, k, 2) \coarser \cdots \coarser \cP(A, k, r) \coarser \cdots \coarser \cP(A, k, \infty) = \cl(\cP(A, k)). 
\end{align}
We write $x \stackrel{P}{\sim} y$ if $x$ and $y$ are in the same part of $P$. 

We define an invariant to record the result of the $k$-WL coloring. 
\begin{align}
    I_{A, k}(t) \coloneqq \sum_{r=0}^\infty t^r M_{A, k}^r,
\end{align}
where
\begin{align}
    M_{A, k}^r = \sum_{(i_1 \ldots i_k) \in V^k} X_{A, k}^r (i_1 \ldots i_k).
\end{align}

Let $A$ and $B$ be adjacency matrices of graphs $G$ and $H$ respectively. 
We say $G \cong_{\text{$k$-WL}} H$ if $I_{A, k}(t) = I_{B, k}(t)$. 
We first characterize when $I_{A, k}(t) = I_{B, k}(t)$.

\begin{proposition}[{Analogue of \cite[Proposition 4]{alzaga2010SpectraSymmetricPowers}}]\label{prop:breeze}
    Let $A, B : V \times V \to C$ be two square matrices. 
    Then $I_{A, k}(t) = I_{B, k}(t)$ if and only if there exists a permutation $\sigma$ of $k$-tuples in $V^k$ such that $X_{A, k}^r (i_1 \ldots i_k) = X_{B, k}^r (\sigma(i_1 \ldots i_k))$ for all $r \geq 1$.
    In particular, 
    \begin{align}
        \tp_A(i_1 \ldots i_k) = \tp_B(\sigma(i_1 \ldots i_k)).
    \end{align}
\end{proposition}

\begin{proof}
    The `\emph{if}' part is obvious. 
    Conversely, assume $I_{A, k}(t) = I_{B, k}(t)$. 
    The coefficient of $t^r$ for $r = \card{V}^k$ implies the existence of a permutation $\sigma$ on the set $V^k$ of $k$-tuples such that
    \begin{align}\label{eq:another}
        X_{A, k}^{\card{V}^k} (i_1 \ldots i_k) = X_{B, k}^{\card{V}^k} (\sigma(i_1 \ldots i_k))
    \end{align}
    Then
    \begin{align}\label{eq:made}
        X_{A, k}^r (i_1 \ldots i_k) = X_{B, k}^r (\sigma(i_1 \ldots i_k))
    \end{align}
    holds for all $1 \leq r \leq \card{V}^k$. 
    Since the WL refinement stabilizes after at most $\card{V}^k$ iterations,~\cref{eq:made} holds for $r \geq \card{V}^k$.
\end{proof}

In fact, if we apply partial refinement during the process of WL algorithm, we will end up with the same result. 

\begin{theorem}\label{thm:black}
    Let $A : V \times V \to C$ be a square matrix and let $P$ be the partition of $V^2$ induced by $\tp_A$. 
    Suppose $P'$ is a partition of $V^2$ induced by $A' : V \times V \to C'$ such that $P \coarser P' \coarser \cl(P)$. 
    Then $\cl(P') = \cl(P)$. 
    In particular, we have $X_{A, 2}^r (i_1, i_2) = X_{A, 2}^r (j_1, j_2)$ for all $r \geq 1$ if and only if $X_{A', 2}^r (i_1, i_2) = X_{A', 2}^r (j_1, j_2)$ for all $r \geq 1$.
\end{theorem}

\begin{proof}
    If one acknowledges that $\cl$ is an idempotent operator, then we have $\cl(P) \coarser \cl(P') \coarser \cl(\cl(P)) = \cl(P)$. 
    Therefore, $\cl(P) = \cl(P')$. 
    Here we give a proof step by step.
    It is clear that $\cl(P) \coarser \cl(P')$. 
    We only need to show that $\cl(P') \coarser \cl(P)$. 
    Suppose $(i_1, i_2) \stackrel{\cl(P)}{\sim} (j_1, j_2)$.
    Then $X_{A, 2}^r (i_1, i_2) = X_{A, 2}^r (j_1, j_2)$ for all $r \geq 1$. 
    By~\cref{prop:breeze}, there exists a permutation of $V^2$ such that
    \begin{align}
        X_{A, 2}^{r+1} (i_1, i_2) = X_{A, 2}^{r+1} (\sigma(i_1, i_2))
    \end{align}
    for all $r \geq 1$, where $(j_1, j_2) = \sigma(i_1, i_2)$.
    We have
    \begin{align}
        X_{A, 2}^{r} (i_1, i_2) = X_{A, 2}^{r} (j_1, j_2),
    \end{align}
    and
    \begin{align}
        \sum_{m \in V} (\tp(i_1, i_2, m), X_{A,2}^r (i_1,m), X_{A,2}^r (m,i_2)) = \sum_{m \in V} (\tp(j_1, j_2, m), X_{B,2}^r (j_1,m), X_{B,2}^r (m,j_2)).
    \end{align}
    Therefore, there exists a permutation $\tau = \tau(\sigma, i_1, i_2)$ such that 
    \begin{align}
        \tp_A(i_1, i_2, i_3) &= \tp_A(j_1, j_2, \tau(i_3)), \\
        X_{A,2}^r (i_1, i_3) &= X_{A, 2}^r (j_1, \tau(i_3)), \label{eq:dance} \\
        X_{A, 2}^r (i_3, i_2) &= X_{A, 2}^r (\tau(i_3), j_2), \label{eq:they}
    \end{align}
    for all $i_3 \in V$. 
    Since $\tp_A(i_1, i_2, i_3) = \tp_A(j_1, j_2, \tau(i_3))$, we know that the first entry ($0$ or $1$) in $T(A', 3)_{p, q}$ is identical for $(i_1, i_2, i_3)$ and $(j_1, j_2, \tau(i_3))$. 
    Now we consider the second entry. 
    Note that we have $X_{A, 2}^r (i_1, i_2) = X_{A, 2}^r (j_1, j_2)$, $X_{A, 2}^r (i_1, i_3) = X_{A, 2}^r (j_1, \tau(i_3))$, and $X_{A, 2}^r (i_3, i_2) = X_{A, 2}^r (\tau(i_3), j_2)$ for all $r \geq 1$. 
    Take $r$ sufficiently large, and we obtain $(i_1, i_2) \stackrel{\cl(P)}{\sim} (j_1, j_2)$, $(i_1, i_3) \stackrel{\cl(P)}{\sim} (j_1, \tau(i_3))$, and $(i_3, i_2) \stackrel{\cl(P)}{\sim} (\tau(i_3), j_2)$. 
    Hence, $(i_1, i_2) \stackrel{P'}{\sim} (j_1, j_2)$, $(i_1, i_3) \stackrel{P'}{\sim} (j_1, \tau(i_3))$, and $(i_3, i_2) \stackrel{P'}{\sim} (\tau(i_3), j_2)$. 
    Moreover, take $i_2 = i_1$ in~\cref{eq:dance,eq:they}, and we conclude that $X_{A, 2}^r (i_1, i_3) = X_{A, 2}^r (j_1, j_3)$ if and only if $X_{A, 2}^r (i_3, i_1) = X_{A, 2}^r (j_3, j_1)$. 
    Therefore, $\tp_{A'} (i_1, i_2, i_3) = \tp_{A'} (j_1, j_2, \tau(i_3))$. 
    So $X_{A', 2}^1 (i_1, i_2) = X_{A', 2}^1 (j_1, j_2) = X_{A', 2}^1 (\sigma(i_1, i_2))$. 
    Next, we prove that $X_{A', 2}^r (i_1, i_2) = X_{A', 2}^r (\sigma(i_1, i_2))$. 
    This is indeed the case since
    \begin{align}
        \tp_{A'}(i_1, i_2, i_3) &= \tp_{A'}(j_1, j_2, \tau(i_3)), \\
        X_{A', 2}^r (i_1, i_3) &= X_{A', 2}^r (j_1, \tau(i_3)),  \\
        X_{A', 2}^r (i_3, i_2) &= X_{A', 2}^r (\tau(i_3), j_2), 
    \end{align}
    holds by induction. 
    In other words, $(i_1, i_2) \stackrel{\cl(P')}{\sim} \sigma(j_1, j_2)$. 
    We conclude that $\cl(P) = \cl(P')$. 
\end{proof}

The $2$-WL algorithm is stronger than distinguishing by matrix powers. 

\begin{theorem}[{Analogue of \cite[Theorem 3]{alzaga2010SpectraSymmetricPowers}}]\label{thm:mainly}
    Let $A, B : V \times V \to C$ be two square matrices, where $C$ is a (not necessarily commutative) ring. 
    If $X_{A, 2}^r (i,j) = X_{B, 2}^r (p,q)$, then $(A^r)(i,j) = (B^r)(p, q)$. 
\end{theorem}

\begin{proof}
    We prove by induction on the number of iterations. 
    The case $r=1$ is trivial. 
    Suppose the claim is true for $r$. 
    Then for $r+1$ we have
    \begin{align}
        X_{A, 2}^{r+1} (i,j) = X_{B, 2}^{r+1} (p,q)
    \end{align}
    By~\cref{eq:iter_rplus1}, we have
    \begin{align}
        \sum_{m \in V} (\tp(i, j, m), X_{A,2}^r (i,m), X_{A,2}^r (m,j)) = \sum_{m \in V} (\tp(p, q, m), X_{B,2}^r (p,m), X_{B,2}^r (m,q))
    \end{align}
    Therefore there exists a permutation $\tau$ on $V$ such that
    \begin{align}
        \tp_A(i, j, m) &= \tp_B(p, q, \tau(m)) \\
        X_{A,2}^r (i, m) &= X_{B, 2}^r (p, \tau(m)) \\
        X_{B, 2}^r (m, j) &= X_{B, 2}^r (\tau(m), q)
    \end{align}
    By the induction hypothesis, we have
    \begin{align}
        A(i, m) &= B(p, \tau(m)) \\
        A(m, j) &= B(\tau(m), q) \\
        (A^r) (i, m) &= (B^r)(p, \tau(m)) \\
        (A^r) (m, j) &= (B^r)(\tau(m), q)
    \end{align}
    By summing over $m$, we get
    \begin{align}
        \sum_m A(i, m) (A^r)(m, j) = \sum_m B(p, m) (B^r)(m, q)
    \end{align}
    In other words, $(A^{r+1})(i,j) = (B^{r+1})(p,q)$. 
\end{proof}

Next, we show that if $G \cong_{\text{$2$-WL}} H$, then their multivariate graph spectra are identical.

\begin{theorem}[{Analogue of~\cite[Theorem 4]{alzaga2010SpectraSymmetricPowers}}]\label{thm:invtospec}
    Let $A_1$ and $B_1$ be two real symmetric matrices of order $n$. 
    Let $e_i$, $i = 1,2, \ldots, p$ be zero-one vectors of length $n$ such that the positions of the ones are disjoint. 
    Let $J_{i,j} = e_i e_j^\top$ for $i,j = 1, 2, \ldots, p$ and $k = 1 + p^2$.
    Let $\bm{A} = (A_1, A_2, \ldots, A_k) = (A_1, J_{1,1}, J_{1,2}, \ldots, J_{p,p})$ and $\bm{B} = (B_1, B_2, \ldots, B_k) = (B_1, J_{1,1}, J_{1,2}, \ldots, J_{p,p})$. 
    If $I_{W(\bm{A}, \bm{s}), 2}(t) = I_{W(\bm{B}, \bm{s}), 2}(t)$, then $\phi_{\bm{A}}(\bm{s}; t) = \phi_{\bm{B}}(\bm{s}; t)$.
\end{theorem}

\begin{proof}
    Assume $I_{W(\bm{A}, \bm{s}), 2}(t) = I_{W(\bm{B}, \bm{s}), 2}(t)$. 
    By~\cref{prop:breeze}, there exists a permutation $\sigma$ of the set $V^2$ of $2$-tuples such that for every $2$-tuple $(i,j)$, it holds
    \begin{align}
        X_{W(\bm{A}, \bm{s}), 2}^r(i,j) = X_{W(\bm{B}, \bm{s}), 2}^r (\sigma(i,j))
    \end{align}
    for all $r \geq 1$. 
    For $r = 1$, we have
    \begin{align}
        \tp_{W(\bm{A}, \bm{s})}(i, j) = \tp_{W(\bm{B}, \bm{s})}(\sigma(i,j)).
    \end{align}
    If $i = j$, then $\sigma$ sends the diagonal of $X_{W(\bm{A}, \bm{s}), 2}^r$ to the diagonal of $X_{W(\bm{B}, \bm{s}), 2}$, namely
    \begin{align}
        \sigma(i, i) = (p, p)
    \end{align}
    for some $p \in V$. 
    Therefore, 
    \begin{align}
        \sum_{i} X_{W(\bm{A}, \bm{s}), 2}^r (i, i) = \sum_{i} X_{W(\bm{B}, \bm{s}), 2}^r (\sigma(i, i)).
    \end{align}
    By~\cref{thm:mainly}, we get
    \begin{align}
        \sum_{i} (W(\bm{A}, \bm{s}))^r (i, i) = \sum_{i} (W(\bm{B}, \bm{s}))^r (\sigma(i, i)).
    \end{align}
    In other words $\Tr(W(\bm{A}, \bm{s})^r) = \Tr((W(\bm{B}, \bm{s}))^r)$ for all $r \geq 1$. 
    Then $\phi_{\bm{A}}(\bm{s}; t) = \phi_{\bm{B}}(\bm{s}; t)$.
\end{proof}

We consider several partial refinements during the process of the $2$-WL. 

\begin{theorem}\label{thm:Aprimes}
    Let $G$ and $H$ be two graphs on $V$ sharing the same degree sequences. 
    Without loss of generality, suppose $\deg(v) = \deg_{G}(v) = \deg_{H}(v)$ for all $v \in V$. 
    Therefore, we have the degree decomposition of vertices $V = \bigsqcup_{i=1}^p V_i$.
    Let $A$ and $B$ be the adjacency matrices of $G$ and $H$ respectively. 
    Let $e_i \in \RR^V$, $i = 1,2, \ldots, p$ be zero-one vectors such that $e_i(v) = 1$ if and only if $v \in V_i$. 
    Let $J_{i,j} = e_i e_j^\top$ for $i,j = 1, 2, \ldots, p$.
    Let $D_i(u, v) = \begin{cases}
        1, & u = v \in V_i, \\
        0, & \textotherwise.
    \end{cases}$
    \begin{enumerate}
        \item Define $A' = s_0 A + s_1 D_1 + \cdots + s_p D_p$ and $B' = s_0 B + s_1 D_1 + \cdots + s_p D_p$. 
        \item Define $A'' = s_0 A_0 + s_1 J_{1,1} + \cdots + s_p J_{p, p}$ and $B'' = s_0 B + s_1 J_{1,1} + \cdots + s_p J_{p, p}$.  
        \item Define $A''' = s_0 A_0 + s_{1,1} J_{1,1} + s_{1,2} J_{1, 2} + \cdots + s_{p,p} J_{p, p}$ and $B''' = s_0 B + s_{1,1} J_{1,1} + s_{1,2} J_{1, 2} + \cdots + s_{p,p} J_{p, p}$.
    \end{enumerate}
     
    Then the following are equivalent. 
    \begin{enumerate}
        \item $I_{A, 2}(t) = I_{B, 2}(t)$.
        \item $I_{A', 2}(t) = I_{B', 2}(t)$.
        \item $I_{A'', 2}(t) = I_{B'', 2}(t)$.
        \item $I_{A''', 2}(t) = I_{B''', 2}(t)$.
    \end{enumerate}
\end{theorem}

\begin{proof}
    Let $P, P', P'', P'''$ be the partition induced by $A, A', A'', A'''$ respectively. 
    By~\cref{thm:black}, we only need to show that $P \coarser P' \coarser \cl(P)$ and $P \coarser P'' \coarser P''' \coarser \cl(P)$ since the case for $B$ is identical. 
    Suppose $(u, v) \stackrel{\cl(P)}{\sim} (z,w)$.
    By~\cref{prop:breeze}, there exists a permutation $\sigma$ of $2$-tuples such that
    \begin{align}
        X_{A, 2}^r (u, v) = X_{A, 2}^r(z, w)
    \end{align}
    for all $r \geq 1$, where $(z, w) = \sigma(u, v)$ for $(u, v) \in V^2$.
    By~\cref{thm:mainly}, we have $(A^r)(u, v) = (A^r)(z, w)$. 
    Note that $\sigma$ maps the diagonal of $X_{A, 2}^r$ to the diagonal of $X_{A, 2}^r$.
    In particular, we have $\deg(u) = (A^2)(u, u) = (A^2)(z, z) = \deg(z)$. 
    Therefore, we have $\sigma (\diag(V_i^2)) = \diag(V_i^2)$ for $i = 1,2, \ldots, p$, where $\diag(V_i^2) \coloneqq \set{(v,v) \in V_i^2}$. 
    We conclude that $P \coarser P' \coarser \cl(P)$. 

    It is clear that $\sigma(\offdiag(V^2)) = \offdiag(V^2)$, where $\offdiag(V^2) \coloneqq \set{(u,v) \in V^2 \mid u \neq v}$. 
    We claim that $\sigma(\offdiag(V_i^2)) = \offdiag(V_i^2)$ for $i = 1,2, \ldots, p$, where $\offdiag(V_i^2) \coloneqq \set{(u,v) \in V_i^2 \mid u \neq v}$. 
    Suppose $(u_1, u_2) \in \offdiag(V_i^2)$ and $\sigma(u_1, u_2) = (z_1, z_2)$. 
    Take $r = 3$, and expand
    \begin{align}
        X_{A, 2}^3 (u, v) = X_{A, 2}^3 (z, w).
    \end{align}
    We get
    \begin{align}
        & \sum_{m \in V} (\tp_{A}(u_1, u_2, m), X_{A, 2}^2 (u_1, m), X_{A,2}^2 (m, u_2)) \\
        =& \sum_{m \in V} (\tp_{A}(z_1, z_2, m), X_{A, 2}^2 (z_1, m), X_{A,2}^2 (m, z_2)).
    \end{align}
    Therefore, there exists a permutation $\tau$ of $V$ such that
    \begin{align}
        & (\tp_{A}(u_1, u_2, m), X_{A, 2}^2 (u_1, m), X_{A,2}^2 (m, u_2)) \\
        =& (\tp_{A}(z_1, z_2, \tau(m)), X_{A, 2}^2 (z_1, \tau(m)), X_{A,2}^2 (\tau(m), z_2)).        
    \end{align}
    We claim $\tau(u_1) = z_1$ and $\tau(u_2) = z_2$. 
    In fact, $X_{A, 2}^2 (u_1, u_1) = X_{A, 2}^2 (z_1, \tau(u_1))$ implies $\tp_{A}(u_1, u_1) = \tp_{A}(z_1, \tau(u_1))$. 
    Hence, $\tau(u_1) = z_1$. 
    Similarly, $\tau(u_2) = z_2$. 
    In particular, we get $X_{A, 2}^2 (u_1, u_1) = X_{A, 2}^2 (z_1, z_1)$, thus $\deg(u_1) = A^2 (u_1, u_1) = A^2 (z_1, z_1) = \deg(z_1)$. 
    Therefore, $z_1 \in V_i$. 
    Similarly, $z_2 \in V_i$. 
    We conclude that $\sigma(\offdiag(V_i^2)) = \offdiag(V_i^2)$, thus $\sigma(V_i^2) = V_i^2$. 

    Next, we claim that $\sigma(V_i \times V_j) = V_i \times V_j$ for $1 \leq i \neq j \leq p$. 
    Suppose $(u_1, u_2) \in V_i \times V_j$ and $\sigma(u_1, u_2) = (z_1, z_2)$. 
    By a similar argument we get $z_1 \in V_i$ and $z_2 \in V_j$. 
    Therefore, $\sigma(V_i \times V_j) = V_i \times V_j$ for $1 \leq i \neq j \leq p$. 
    Hence, we obtain $P \coarser P'' \coarser P''' \coarser \cl(P)$.
\end{proof}

\section{Proof of main theorems}\label{sec:proofs}

Now we are prepared to prove the main theorems. 

\begin{proof}[{Proof of~\cref{thm:2WLtoGBLS}}]
    Suppose $G \cong_{\text{2-WL}} H$. 
    Then $I_{A, 2}(t) = I_{B, 2}(t)$. 
    By~\cref{thm:Aprimes}, we have $I_{A', 2}(t) = I_{B', 2}(t)$, $I_{A'', 2}(t) = I_{B'', 2}(t)$, $I_{A''', 2}(t) = I_{B''', 2}(t)$. 
    By~\cref{thm:invtospec}, we have $\phi_{A'}(s, t) = \phi_{B'}(s, t)$ (namely $G \cong_{GDLS} H$), $\phi_{A''}(s, t) = \phi_{B''}(s, t)$ (namely $G \cong_{GBDLS} H$), and $\phi_{A'''}(s, t) = \phi_{B'''}(s, t)$ (namely $G \cong_{GBLS} H$). 
\end{proof}

\begin{proof}[{Proof of~\cref{thm:GBLStodn}}]
    Suppose $G \cong_{GBLS} H$. 
    By~\cref{thm:GBLS} we have $Q^\top \widetilde{W}_{A} = \widetilde{W}_B$. 
    Let $\widetilde{W}_A = U \Sigma V$ be the Smith decomposition of $\widetilde{W}_A$. 
    Then $Q^\top U \Sigma = \widetilde{W}_B V^\top$. 
    By~\cref{lem:divide_dn}, we have $\ell(Q^\top U) \mid d_n(\Sigma) = d_n(\widetilde{W}_A)$. 
    We conclude that $\ell(Q) = \ell(Q^\top) = \ell(Q^\top U) \mid d_n(\widetilde{W}_A)$ by~\cref{lem:multiply_unimodular}. 
\end{proof}

\section{Discussion}

We test the power of~\cref{thm:main} by computing 10000 random graphs $G \in \mathbb{G}(n, 1/2)$ of order $n$ for $n = 10, 20, 30, 40, 50$. 
The result of one experiment is summarized in~\cref{tbl:softly}.

\begin{table}[htbp]
    \centering
    \begin{tabular}{cccc}
        \hline
        $n$ & $\#\set{G : d_n(\widetilde{W}_A) = 0}$ & $\#\set{G : d_n(\widetilde{W}_A) = 1}$ & $\#\set{G : \text{~\cref{thm:main} is applicable}}$\\
        \hline
        10 & 2424 & 6473 & 6510 \\
        20 & 9 & 9511 & 9517 \\
        30 & 0 & 9910 & 9910 \\
        40 & 0 & 9979 & 9979 \\
        50 & 0 & 9989 & 9989 \\
        \hline
    \end{tabular}
    \caption{Test the power of~\cref{thm:main}}\label{tbl:softly}
\end{table}

The generalized block Laplacian spectrum adopts $k = 1 + p^2$ variables, which could be very large. 
We would like to minimize the number of variables.

In fact, even if we consider a sparser partition, the result of the experiments is still satisfactory. 
Let $\delta_1 > \delta_2 > \cdots > \delta_p$ be the distinct degrees in a graph $G$. 
Let $r = \ceil{\log_2 n} + 1$. 
We consider the partition $V = \bigsqcup_{i=1}^r V_i$, where $V_i$ is the set of vertices of degree $\delta_i$ for $i = 1,2, \ldots, r - 1$ and $V_r = V \setminus \bigsqcup_{i=1}^{r-1} V_i$ the set of rest vertices (if $p > r$ then we fill by empty subset). 
Then with high probability, we have $|V_1| = |V_2| = \cdots = |V_{r-1}| = 1$~\cite[Theorem 3.6]{babai1980RandomGraphIsomorphism}.
Let $e_i$ be the characteristic vector of $V_i$ and $J_{i,j} = e_i e_j^\top$ for $1 \leq i,j \leq r$. 
We call the multivariate graph spectrum for $\bm{A} = (A_G, J_{1,1}, J_{1,2}, \ldots, J_{r,r})$ the truncated generalized block Laplacian spectrum. 
We consider 
\begin{align}
    \widehat{W}_A \coloneqq [e_1, A e_1, \ldots, A^{n-1} e_1, e_2, A e_2, \ldots, A^{n-1} e_2, \ldots, e_r, A e_r, \ldots, A^{n-1} e_r].
\end{align}
The analogue of~\cref{thm:GBLStodn} is obvious. 
And the result of one experiment is summarized in~\cref{tbl:cattle}.

\begin{table}[htbp]
    \centering
    \begin{tabular}{cccc}
        \hline
        $n$ & $\#\set{G : d_n(\widehat{W}_A) = 0}$ & $\#\set{G : d_n(\widehat{W}_A) = 1}$\\
        \hline
        10 & 2462 & 6407 \\
        20 & 3 & 9260 \\
        30 & 0 & 9383 \\
        40 & 0 & 9679 \\
        50 & 0 & 9688 \\
        \hline
    \end{tabular}
    \caption{Test the power of truncated generalized block Laplacian spectrum}\label{tbl:cattle}
\end{table}

One may wonder whether we can further reduce the number of variables. 
If we consider the partition of vertices by their degree modulo $4$, then the result of experiments is not satisfactory.

Another direction is to discard the off-diagonal block all-one matrices $J_{i, j}$. 
In this case, the $\RR$-linear combination of matrices are all real symmetric, hence the eigenvalues are real. 
But extra conditions, including but not limited to connectedness, may be necessary for the existence of a fixed orthogonal matrix, as in Theorem~\ref{thm:GBLS}.

We raise several conjectures for future research. 

\begin{conjecture}\label{conj:save}
    Almost all graphs are determined by their truncated generalized block Laplacian spectra. 
\end{conjecture}

\begin{conjecture}\label{conj:rush}
    Almost all graphs are determined by their generalized block Laplacian spectra (GBLS). 
\end{conjecture}

\cref{conj:rush} is a weaker version of Haemers' conjecture.
And~\cref{conj:save} implies~\cref{conj:rush}. 
We do not know how to prove these conjectures. 
A rough idea is to show that 
\begin{align}
    m_i = \det[e_1, e_2, \ldots, e_r, A e_i, \ldots, A^{n-r} e_i]
\end{align}
are almost independent random integers for $1 \leq i \leq r-1$. 
The probability that $k$ random positive integers from $\set{1,2, \ldots, n}$ are coprime is $1/\zeta(k) + \mathcal{O}(1/n)$, where $\zeta(k)$ is the Riemann zeta function~\cite{nymann1972ProbabilityThatPositive}. 
Note that $1/\zeta(\ceil{\log_2 50}) \approx 0.9830$.

At last, we would like to point out that though we introduce the multivariate graph spectrum in this paper, it in fact can be reduced to the spectrum of a fixed linear combination of graph matrices. 
The following argument is improved by an anonymous referee. 
Note that integral matrices are useful to avoid the precision issue of floating numbers. 

\begin{lemma}\label{lem:sign}
    Let $f_1, f_2, \ldots, f_m \in \RR[\bm{x}]$ be distinct homogeneous polynomials of degree $d$ in $\bm{x} = (x_1, x_2, \ldots, x_k)$. 
    Then there exists $\bm{0} \neq \bm{c} = (c_1, c_2, \ldots, c_k) \in \ZZ^k$ such that $f_1(\bm{c})$, $f_2(\bm{c})$$, \ldots, f_m(\bm{c})$ are all distinct. 
\end{lemma}

\begin{proof}
    Note that $g = \prod_{1 \leq i < j \leq m}(f_i - f_j) \in \RR[\bm{x}]$ is a nonzero polynomial. 
    Since $\QQ^k$ in dense in $\RR^k$, there exists $\bm{c} = (c_1, c_2, \ldots, c_k) \in \QQ^k$ such that $g(\bm{c}) \neq 0$. 
    Hence, $f_1(\bm{c}), f_2(\bm{c}), \ldots, f_m(\bm{c})$ are all distinct. 
    Let $N$ be a nonzero integer such that $N\bm{c} \in \ZZ^k$.
    Since $f_1, f_2, \ldots, f_m$ are homogeneous polynomials of degree $d$, we have that $f_1(N\bm{c}), f_2(N\bm{c}), \ldots, f_m(N\bm{c})$ are all distinct.
\end{proof}

\begin{theorem}\label{thm:wrapped}
    Let $n, k$ be positive integers. 
    For every graph $G$ with $n$ vertices, let $\bm{A}_G$ be a $k$-tuple of $n \times n$ graph matrices over $\RR$. 
    Then there exists $\bm{s}_0 \in \ZZ^k$ depending only on $n$ such that for every pair of graphs $G, H$ with $n$ vertices $\phi_{\bm{A}_G}(\bm{s}; t) = \phi_{\bm{A}_H}(\bm{s}; t)$ holds if and only if $\phi_{\bm{A}_G}(\bm{s}_0; t) = \phi_{\bm{A}_H}(\bm{s}_0; t)$ holds. 
    Namely, the spectra of $\bm{A}_G(\bm{s}_0)$ and $\bm{A}_H(\bm{s}_0)$ are identical implies that the spectra of $\bm{A}_G(\bm{s})$ and $\bm{A}_H(\bm{s})$ are identical for arbitrary $\bm{s}$. 
\end{theorem}

\begin{proof}
    Consider the set of polynomials $\set{\phi_{\bm{A}_G}(\bm{s}; t) : G \text{ is a graph with } n \text{ vertices}}$. 
    The theorem follows by~\cref{lem:sign}.
\end{proof}

\begin{problem}
    Find an exact $\bm{s}_0 \in \ZZ^k$ in~\cref{thm:wrapped} (for certain tuples of graph matrices).
\end{problem}

\section*{Acknowledgements}

Wei WANG is supported in part by National Key Research and Development Program of China 2023YFA1010203 and National Natural Science Foundation of China (No. 12371357). 
Da ZHAO is supported in part by the National Natural Science Foundation of China (No. 12471324, No. 12501459, No. 12571353), and the Natural Science Foundation of Shanghai, Shanghai Sailing Program (No. 24YF2709000).  
The authors thank Qing XIANG for helpful discussions. 
The authors thank the anonymous referee for the improvement of the paper. 

\bibliographystyle{alpha}
\bibliography{ref}

\end{document}

%% file: relation.tex
\begin{tikzpicture}[
    rel/.style={
      % The shape:
      rectangle,
      % The size:
      minimum size=6mm,
      % The border:
      very thick,
      draw=black,         
                                    % and that mixed with 50% white
      % The filling:
    %   top color=white,              % a shading that is white at the top...
    %   bottom color=red!50!black!20, % and something else at the bottom
      % Font
    %   font=\itshape
    }]

    \node (iso) at (0,0) [rel] {$G \cong_{\text{iso}} H$};
    \node (kWL) at (0, -2) [rel] {$G \cong_{k-\text{WL}} H$};
    \node (k3WL) at (0, -4) [rel] {$G \cong_{3-\text{WL}} H$};
    \node (k2WL) at (0, -6) [rel] {$G \cong_{2-\text{WL}} H$};
    \node (k1WL) at (0, -9) [rel] {$G \cong_{1-\text{WL}} H$};

    \node (GBLS) at (3, -5) [rel] {$G \cong_{GBLS} H$};
    \node (GBDLS) at (3, -6) [rel] {$G \cong_{GBDLS} H$};
    \node (GS) at (3, -7) [rel] {$G \cong_{GS} H$};
    \node (GDLS) at (3, -8) [rel] {$G \cong_{GDLS} H$};

    \node (frac) at (3, -9) [rel] {$G \cong_{\text{frac}} H$};

    \node (qc) at (6, -4) [rel] {$G \cong_{qc} H$};

    \node (S) at (6, -8) [rel] {$G \cong_{S} H$};

    \node (Fall) at (10, 0) [rel] {$\mathrm{hom}(F,G) = \mathrm{hom}(F, H) \forall F$};
    \node (Fplanar) at (10, -4) [rel] {$\mathrm{hom}(P,G) = \mathrm{hom}(P, H) \forall P$};
    \node (Fcycle) at (10, -8) [rel] {$\mathrm{hom}(C,G) = \mathrm{hom}(C, H) \forall C$};
    \node (Ftree) at (10, -9) [rel] {$\mathrm{hom}(T,G) = \mathrm{hom}(T, H) \forall T$};

    \node (Ck++) at (-3, -2) [rel] {$G \cong_{\mathrm{C}^{k+1}} H$};
    \node (C4) at (-3, -4) [rel] {$G \cong_{\mathrm{C}^{4}} H$};
    \node (C3) at (-3, -6) [rel] {$G \cong_{\mathrm{C}^{3}} H$};
    \node (C2) at (-3, -9) [rel] {$G \cong_{\mathrm{C}^{2}} H$};

    \draw[->] (iso) -- (kWL);
    \draw[->] (kWL) -- (k3WL);
    \draw[->] (k3WL) -- (k2WL);
    \draw[->] (k2WL) -- (k1WL);
    \draw[->] (GBLS) -- (GBDLS);
    \draw[->] (GBDLS)-- (GS);

    \draw[->] (k2WL) -- (GBLS);
    \draw[->] (k2WL) -- (GBDLS);
    \draw[->] (k2WL) -- (GS);
    \draw[->] (k2WL) -- (GDLS);

    \draw[->] (GBLS) to[out = 0] (S);
    \draw[->] (GBDLS) to[out = 0] (S);
    \draw[->] (GS) to[out = 0] (S);
    \draw[->] (GDLS) -- (S);

    \draw[->] (Fall) -- (Fplanar);
    \draw[->] (Fplanar) -- (Fcycle);
    \draw[->] (Fplanar) to[out = 0, in = 0] (Ftree);

    \draw[<->] (iso) -- (Fall);
    \draw[<->] (qc) -- (Fplanar);
    \draw[<->] (S) -- (Fcycle);
    \draw[<->] (k1WL) -- (frac);
    \draw[<->] (frac) -- (Ftree);

    \draw[<->] (Ck++) -- (kWL);
    \draw[<->] (C4) -- (k3WL);
    \draw[<->] (C3) -- (k2WL);
    \draw[<->] (C2) -- (k1WL);

    % \draw[->] (GS) to [out = 0, in = 0] (frac);

    \draw[->, dashed] (k3WL) to [out = 0, in = -10] (iso);

    \draw[->, dashed] (GBLS) to [out = 90, in = -10] (iso);
    \draw[->, dashed] (GBDLS) to [out = 20, in = 0] (iso);
    \draw[->, dashed] (GS) to [out = 10, in = 10] (iso);

    \draw[->, dashed] (k1WL) to [out = 150, in = 210] (iso);

\end{tikzpicture}

%% file: ref.bib
@article{wang2016SquarefreeDiscriminantsMatrices,
	title = {Square-free {Discriminants} of {Matrices} and the {Generalized} {Spectral} {Characterizations} of {Graphs}},
	url = {http://arxiv.org/abs/1608.01144},
	urldate = {2021-10-08},
	journal = {arXiv:1608.01144 [math]},
	author = {Wang, Wei and Yu, Tao},
	month = aug,
	year = {2016},
	note = {arXiv: 1608.01144
version: 1},
}

@article{alzaga2010SpectraSymmetricPowers,
  title = {Spectra of Symmetric Powers of Graphs and the {{Weisfeiler}}--{{Lehman}} Refinements},
  author = {Alzaga, Afredo and Iglesias, Rodrigo and Pignol, Ricardo},
  year = {2010},
  month = nov,
  journal = {Journal of Combinatorial Theory, Series B},
  volume = {100},
  number = {6},
  pages = {671--682},
  issn = {0095-8956},
  doi = {10.1016/j.jctb.2010.07.001},
  langid = {english}
}

@article{wang2024HaemersConjectureAlgorithmic,
  title = {Haemers' {{Conjecture}}: {{An Algorithmic Perspective}}},
  shorttitle = {Haemers' {{Conjecture}}},
  author = {Wang, Wei and Wang, Wei},
  year = {2024},
  month = apr,
  journal = {Experimental Mathematics},
  pages = {1--28},
  issn = {1058-6458, 1944-950X},
  doi = {10.1080/10586458.2024.2337229},
  urldate = {2024-10-22},
  langid = {english}
}

@article{cai1992OptimalLowerBound,
  title = {An Optimal Lower Bound on the Number of Variables for Graph Identification},
  author = {Cai, Jin-Yi and F{\"{u}}rer, Martin and Immerman, Neil},
  year = {1992},
  month = dec,
  journal = {Combinatorica},
  volume = {12},
  number = {4},
  pages = {389--410},
  issn = {0209-9683},
  doi = {10.1007/BF01305232},
  langid = {english}
}

@inproceedings{grohe2000IsomorphismTestingEmbeddable,
  title = {Isomorphism Testing for Embeddable Graphs through Definability},
  booktitle = {Proceedings of the Thirty-Second Annual {{ACM}} Symposium on {{Theory}} of Computing},
  author = {Grohe, Martin},
  year = {2000},
  month = may,
  pages = {63--72},
  publisher = {ACM},
  address = {Portland Oregon USA},
  doi = {10.1145/335305.335313},
  urldate = {2024-11-23},
  isbn = {978-1-58113-184-0},
  langid = {english}
}

@article{boppana1987DoesCoNPHave,
  title = {Does Co-{{NP}} Have Short Interactive Proofs?},
  author = {Boppana, Ravi B. and Hastad, Johan and Zachos, Stathis},
  year = {1987},
  month = may,
  journal = {Information Processing Letters},
  volume = {25},
  number = {2},
  pages = {127--132},
  issn = {00200190},
  doi = {10.1016/0020-0190(87)90232-8},
  urldate = {2024-12-12},
  copyright = {https://www.elsevier.com/tdm/userlicense/1.0/},
  langid = {english}
}

@article{schoning1988GraphIsomorphismLow,
  title = {Graph Isomorphism Is in the Low Hierarchy},
  author = {Sch{\"o}ning, Uwe},
  year = {1988},
  month = dec,
  journal = {Journal of Computer and System Sciences},
  volume = {37},
  number = {3},
  pages = {312--323},
  issn = {00220000},
  doi = {10.1016/0022-0000(88)90010-4},
  urldate = {2024-12-12},
  copyright = {https://www.elsevier.com/tdm/userlicense/1.0/},
  langid = {english}
}

@inproceedings{hopcroft1974LinearTimeAlgorithm,
  title = {Linear Time Algorithm for Isomorphism of Planar Graphs (Preliminary Report)},
  booktitle = {Proceedings of the Sixth Annual {{ACM}} Symposium on {{Theory}} of Computing  - {{STOC}} '74},
  author = {Hopcroft, J. E. and Wong, J. K.},
  year = {1974},
  pages = {172--184},
  publisher = {ACM Press},
  address = {Seattle, Washington, United States},
  doi = {10.1145/800119.803896},
  urldate = {2024-12-12},
  langid = {english}
}

@article{tinhofer1991NoteCompactGraphs,
  title = {A Note on Compact Graphs},
  author = {Tinhofer, G.},
  year = {1991},
  month = feb,
  journal = {Discrete Applied Mathematics},
  volume = {30},
  number = {2-3},
  pages = {253--264},
  issn = {0166218X},
  doi = {10.1016/0166-218X(91)90049-3},
  urldate = {2024-12-12},
  copyright = {https://www.elsevier.com/tdm/userlicense/1.0/},
  langid = {english}
}

@article{lovasz1967OperationsStructures,
  title = {Operations with Structures},
  author = {Lov{\'a}sz, L.},
  year = {1967},
  month = sep,
  journal = {Acta Mathematica Academiae Scientiarum Hungaricae},
  volume = {18},
  number = {3-4},
  pages = {321--328},
  issn = {0001-5954, 1588-2632},
  doi = {10.1007/BF02280291},
  urldate = {2024-12-12},
  copyright = {http://www.springer.com/tdm},
  langid = {english}
}

@article{babai1980RandomGraphIsomorphism,
  title = {Random Graph Isomorphism},
  author = {Babai, L{\'a}szl{\'o} and Erd{\"o}s, Paul and Selkow, Stanley M.},
  year = {1980},
  month = aug,
  journal = {SIAM Journal on Computing},
  volume = {9},
  number = {3},
  pages = {628--635},
  issn = {0097-5397, 1095-7111},
  doi = {10.1137/0209047},
  urldate = {2024-10-28},
  langid = {english}
}

@inproceedings{babai1979CanonicalLabellingGraphs,
  title = {Canonical Labelling of Graphs in Linear Average Time},
  booktitle = {20th {{Annual Symposium}} on {{Foundations}} of {{Computer Science}} (Sfcs 1979)},
  author = {Babai, Laszlo and Kucera, Ludik},
  year = {1979},
  month = oct,
  pages = {39--46},
  issn = {0272-5428},
  doi = {10.1109/SFCS.1979.8},
  urldate = {2024-12-12},
  langid = {english}
}

@article{mancinskaQuantumIsomorphismEquivalent2019,
  title = {Quantum Isomorphism Is Equivalent to Equality of Homomorphism Counts from Planar Graphs},
  author = {Man{\v c}inska, Laura and Roberson, David E.},
  year = {2019},
  month = oct,
  journal = {arXiv:1910.06958 [quant-ph]},
  eprint = {1910.06958},
  primaryclass = {quant-ph},
  urldate = {2021-03-25},
  archiveprefix = {arXiv}
}

@article{qiu2023GeneralizedSpectralCharacterizations,
  title = {Generalized Spectral Characterizations of Regular Graphs Based on Graph-Vectors},
  author = {Qiu, Lihong and Ji, Yizhe and Mao, Lihuan and Wang, Wei},
  year = {2023},
  month = apr,
  journal = {Linear Algebra and its Applications},
  volume = {663},
  pages = {116--141},
  issn = {0024-3795},
  doi = {10.1016/j.laa.2023.01.006},
  urldate = {2024-11-01}
}

@article{grohe2025IterationNumberWeisfeilerleman,
	title = {The iteration number of the weisfeiler-leman algorithm},
	volume = {26},
	issn = {1529-3785},
	url = {https://dl.acm.org/doi/10.1145/3708891},
	doi = {10.1145/3708891},
	language = {en},
	number = {1},
	urldate = {2025-08-31},
	journal = {ACM Trans. Comput. Logic},
	author = {Grohe, Martin and Lichter, Moritz and Neuen, Daniel},
	year = {2025},
	pages = {6:1--6:31},
}

@article{johnson1980NoteCospectralGraphs,
	title = {A note on cospectral graphs},
	volume = {28},
	issn = {0095-8956},
	url = {https://www.sciencedirect.com/science/article/pii/0095895680900581},
	doi = {10.1016/0095-8956(80)90058-1},
	language = {en},
	number = {1},
	urldate = {2021-09-08},
	journal = {Journal of Combinatorial Theory, Series B},
	author = {Johnson, Charles R and Newman, Morris},
	month = feb,
	year = {1980},
	pages = {96--103},
}

@article{nymann1972ProbabilityThatPositive,
	title = {On the probability that \textit{k} positive integers are relatively prime},
	volume = {4},
	issn = {0022-314X},
	url = {https://www.sciencedirect.com/science/article/pii/0022314X72900388},
	doi = {10.1016/0022-314X(72)90038-8},
	language = {en},
	number = {5},
	urldate = {2025-09-14},
	journal = {Journal of Number Theory},
	author = {Nymann, J. E},
	month = oct,
	year = {1972},
	pages = {469--473},
}

@misc{immerman2019,
	title = {The \$k\$-dimensional weisfeiler-leman algorithm},
	url = {http://arxiv.org/abs/1907.09582},
	doi = {10.48550/arXiv.1907.09582},
	language = {en},
	urldate = {2025-10-28},
	publisher = {arXiv},
	author = {Immerman, Neil and Sengupta, Rik},
	month = jul,
	year = {2019},
	note = {arXiv:1907.09582 [cs]},
}

@inproceedings{grohe1998FixedpointLogicsPlanar,
	title = {Fixed-point logics on planar graphs},
	url = {https://ieeexplore.ieee.org/document/705639/},
	doi = {10.1109/LICS.1998.705639},
	language = {en},
	urldate = {2025-10-28},
	booktitle = {Proceedings. {Thirteenth} {Annual} {IEEE} {Symposium} on {Logic} in {Computer} {Science} (cat. {No}.98cb36226)},
	author = {Grohe, M.},
	month = jun,
	year = {1998},
	note = {ISSN: 1043-6871},
	pages = {6--15},
}
